\documentclass[twoside,12pt]{article}
\input{amssymb.sty} \catcode`\@=11
\def\numberbysection{\@addtoreset{equation}{section}
        \def\theequation{\thesection.\arabic{equation}}}
\numberbysection

\newcommand{\be}{\begin{equation}}
\newcommand{\ee}{\end{equation}}
\newcommand{\bea}{\begin{eqnarray}}
\newcommand{\eea}{\end{eqnarray}}

\begin{document}

\def\Z{{\mathbb Z}}
\def\N{{\mathbb N}} 
\def\M{{\mathbb M}}
\def\R{{\mathbb R}}
\def\H{{\mathbb H}}
\def\Q{{\mathbb Q}}
\def\C{{\mathbb C}}
\def\i{{\rm i}}


\vskip 1truecm

\centerline{{\huge The algebraic meaning of genus-zero}\footnote{A
contribution to the Moonshine Conference at ICMS, Edinburgh, July 2004.}}

\vskip 2truecm
\centerline{{\large Terry Gannon}}\smallskip
\centerline{{\it Department of Mathematical Sciences,}}
\centerline{{\it University of Alberta}}
\centerline{{\it Edmonton, Canada, T6G 2G1}}
\centerline{{\tt tgannon@math.ualberta.ca}}

\vskip 1.5truecm

\centerline{{\bf Abstract}}\vskip .3truecm

The Conway--Norton conjectures unexpectedly related the Monster with
certain special modular functions (Hauptmoduls).  
Their proof  by Borcherds {\it et al}
was remarkable for demonstrating the rich mathematics implicit there.
Unfortunately Moonshine remained almost as mysterious after the proof as before.
In particular,  a computer check --- as opposed to a general conceptual
argument --- was used to verify the Monster functions equal the 
appropriate modular functions. This, the so-called `conceptual gap',
was eventually filled; we review the solution here. We conclude by speculating
on the shape of a new proof of the Moonshine conjectures.

\vskip 1.5truecm

\section{The conceptual gap}

The main Conway--Norton conjecture \cite{CN} says:

\medskip\noindent{{\bf Theorem 1.}} {\it There is an infinite-dimensional
graded representation $V=\oplus_{n=-1}^\infty V_n$ of the Monster $\M$, such
that the McKay--Thompson series
\be
T_g(\tau):={\rm Tr}_{V}
g\, q^{L_0-1}=\sum_{n\ge -1}c_n(g)\,q^n
\ee
equals the Hauptmodul $J_g$ for
some discrete subgroup $\Gamma_g$ of SL$_2(\R)$.}\medskip

Moreover, each coefficient $c_n(g)$ lies in $\Z$, and $\Gamma_g$ contains the congruence
subgroup $\Gamma_0(N)$ as a normal subgroup, where $N=h\,o(g)$ for some
$h$ dividing gcd$(24,o(g))$ ($o(g)$ is the order of $g\in\M$).
In his ICM talk \cite{B2}, Borcherds outlined the proof of Theorem 1:

\begin{itemize}

\item[(i)] Construction of the Frenkel--Lepowsky--Meurman Moonshine module
$V^\natural$, which is to equal the space $V$;

\item[(ii)] Derivation of recursions for the McKay--Thompson
coefficients $c_n(g)$, such as
\be
c_{4n+2}(g)=c_{2k+2}(g^2)+\sum_{j=1}^k c_j(g^2)\,c_{2k+1-j}(g^2)\qquad
\forall k\ge 1\ ;
\ee

\item[(iii)] From these recursions, prove $T_g=J_g$.
\end{itemize}

The original treatment of step (i) is \cite{FLM}, and is reviewed elsewhere in these
proceedings. Borcherds derived (ii) by first constructing a Lie algebra
out of $V^\natural$, and then computing its twisted denominator identities 
\cite{B1}. It was already
known that Hauptmoduls automatically satisfied these recursions, and that any function
obeying all those recursions was uniquely determined by its first few
coefficients. Thus establishing (iii) merely requires comparing finitely
many coefficients of each $T_g$ with $J_g$ --- in fact, comparing 5 coefficients
for each of the 171 functions suffices \cite{B1}. In this way Borcherds accomplished
(iii) and with it completed the proof of Theorem 1.
The proof successfully established the mathematical richness of the subject, 
and for
his work Borcherds deservedly received a Fields medal (math's highest honour)
in 1998.

This quick sketch hides the technical sophistication of the proof of (i) and (ii).
More recently,
the construction of (i) has been simplified in \cite{Mi}, and the derivations of
the recursions have been simplified in \cite{JLW,KK}. But the biggest weakness
of the proof is hidden in the nearly trivial argument of step (iii).

At the risk of sending shivers down Bourbaki's collective spine, the point
of mathematics is surely not acquiring proofs (just as the point of theoretical physics
is not careful calculations, and that of painting is not the creation of
realistic scenes on canvas). The point of mathematics, like that of any
intellectual discipline, is to find qualitative truths, to abstract out
patterns from the inundation of seemingly disconnected facts. An example is the
algebraic notion of group. Another,
dear to many of us, is the A-D-E meta-pattern: many different classifications
(e.g.\ finite subgroups of SU$_2$, subfactors of small index, the simplest
conformal field theories) fall unexpectedly into the same pattern.
The conceptual explanation for the ubiquity of this meta-pattern --- that is,
the combinatorial fact underlying its various manifestations --- presumably
involves the graphs with largest eigenvalue $|\lambda|\le 2$. 

Likewise, the real challenge of
Monstrous Moonshine wasn't to prove Theorem 1, but rather to understand
what the Monster has to do with modularity and genus-0. The first proof
was due to Atkin, Fong and Smith \cite{Sm}, who by studying the first
100 coefficients of the $T_g$ verified (without constructing it) that there 
existed a (possibly virtual) representation
$V$ of $\M$ obeying Theorem 1.
Their proof is forgotten because it didn't explain anything. 

By contrast, the proof of Theorem 1 by Borcherds {\it et al} is clearly
superior: it explicitly constructs $V=V^\natural$, and
emphasises the remarkable mathematical richness saturating the problem. 
On the other hand, it also fails to 
explain modularity and the Hauptmodul property. The problem
is step (iii): precisely at the point where we want to identify the algebraically
defined $T_g$'s with the topologically defined $J_g$'s, a conceptually empty
computer check of a few hundred coefficients is done. This is called the
conceptual gap of Monstrous Moonshine, and it has an analogue in Borcherds'
proof of Modular Moonshine \cite{B3} and in H\"ohn's proof of `generalised Moonshine' for the Baby
Monster \cite{Ho}. Clearly preferable would be to replace the numerical check 
of \cite{B1} with a more general theorem. 

Next section we review the standard definition of Hauptmodul. In Section 3 we
describe the solution \cite{CG} to the conceptual gap: it
replaces that topological definition of Hauptmodul with an algebraic
one. We conclude the paper with some speculations. Even with the improvements
\cite{JLW,KK,Mi} and especially \cite{CG} to the original proof of Theorem 1, 
the resulting argument still does a poor job explaining Monstrous 
Moonshine.
Moonshine remains mysterious to this day. There is a lot left to do --- for example
establishing Norton's generalised Moonshine \cite{N}, or finding the Moonshine 
manifold \cite{HBJ}. But the greatest task for Moonshiners is to find a second
independent proof of Theorem 1. It would (hopefully) clarify some things
that the original proof leaves murky. In particular, we still don't know
what really is so important about the Monster, that it has such a rich
genus-0 moonshine. To what extent does Monstrous Moonshine determine the
Monster? We turn to this open problem in Section 4.

\section{The topological meaning of the Hauptmodul property}

Just as a {\it periodic} function is a function on a compact {\it real}
curve (i.e.\ on a circle), a {\it modular} function is a function on a
compact {\it complex} curve. More precisely,
let $\Sigma$ be a compact surface. We can regard this as a complex
curve, and thus put on it a complex analytic structure. Up to (biholomorphic
or conformal)
equivalence, there is a unique genus 0 surface (which we can take to be
 the Riemann sphere $\C\cup\{\infty\}$), but there is a
continuum (moduli space) of inequivalent complex analytic structures which can be placed
on a torus (genus 1), a double-torus, etc. For example, this moduli space
for the torus can be naturally identified with $\H/{\rm SL}_2(\Z)$, where
$\H$ is the upper half-plane $\{\tau\in\C\,|\,{\rm Im}(\tau)>0\}$: 
any torus is equivalent
to one of the form $\C/(\Z+\Z\tau)$ for some $\tau\in\H$; 
the tori corresponding to $\tau$ and $\tau'$ are themselves equivalent
iff $\tau'=(a\tau+b)/(c\tau+d)$ for some $\left(\matrix{a&b\cr
c&d}\right)\in{\rm SL}_2(\Z)$.

We learn from both geometry and physics that we should study a space through the 
functions (fields) that live on it, which respect the relevant
properties of the space. Therefore we should consider meromorphic
functions $f:\Sigma\rightarrow \C$
(a meromorphic function is holomorphic everywhere, except for isolated
finite poles; we would have preferred $f$ to be holomorphic, but then
$f$ would be constant). These $f$ are called {\it modular functions}, and as
they are  as important as complex curves, they should be central
to mathematics. 

We know all about meromorphic functions for $\C$: these include
rational functions (i.e.\ quotients of
polynomials), together with transcendental functions such as $\exp(z)$ and
$\cos(z)$. But both  $\exp(z)$ and $\cos(z)$ have an essential
singularity at $\infty$. In fact, the modular functions for the Riemann
sphere are the rational functions in $z$. By contrast, the modular functions
for the other compact surfaces will be rational functions in two generators,
where those two generators satisfy a polynomial relation. For example, the
modular functions for a torus are generated by the Weierstrass ${\mathfrak p}$-function
and its derivative, and ${\mathfrak p}$ and ${\mathfrak p}'$ satisfy the cubic 
equation defining the torus.
In this way, the sphere is distinguished from all other compact surfaces.

There are three possible geometries in two dimensions: Euclidean, spherical
and hyperbolic. The most important of these, in any sense, is
the hyperbolic one. The upper half-plane is a model for it: its `lines'
consist of vertical half-lines and semi-circles, its infinitesimal metric
is d$s=|{\rm d}\tau|/{\rm Im}(z)$, etc. As with Euclidean geometry, `lines'
are the paths of shortest distance. Just as the Euclidean plane $\R^2$ has
a circular horizon (one infinite point for every angle $\theta$), so does the
hyperbolic plane $\H$, and it can be identified with $\R\cup\{\i\infty\}$.
 The group of isometries (geometry-preserving
transformations $\H\rightarrow\H$) is SL$_2(\R)$, which acts on $\H$ as
fractional linear transformations $\tau\mapsto{a\tau+b\over c\tau+d}$, and
sends the horizon to itself.

It turns out that any compact surface $\Sigma$ can be realised (in infinitely
many different ways) as the compactification of the space  $\H/\Gamma$ of orbits,
for some discrete subgroup $\Gamma$ of SL$_2(\R)$. The compactification
amounts to including finitely many $\Gamma$-orbits of horizon points.
The most important example is $\H/{\rm SL}_2(\Z)$, which can be identified with
the sphere with one puncture, the puncture corresponding to the single compactification
orbit $\Q\cup\{\i\infty\}$. Those points $\Q\cup\{\i\infty\}$ are called
cusps.
The groups $\Gamma$ of greatest interest in number theory, and to us, are those
commensurable to SL$_2(\Z)$: i.e.\ $\Gamma\cap{\rm SL}_2(\Z)$ is also an
infinite discrete group with finite index in both $\Gamma$ and SL$_2(\Z)$.
Their compactification points will again be the cusps $\Q\cup\{\i\infty\}$.
Examples of such groups are SL$_2(\Z)$ itself, as well as its subgroups
\bea
\Gamma(N)&=&\,\{A\in{\rm SL}_2(\Z)\,|\, A\equiv\pm I\ ({\rm mod}\
N)\}\ ,\\
\Gamma_0(N)&=&\,\left\{\left(\matrix{a&b\cr c&d}\right)\in{\rm SL}_2(\Z)\,\bigm|\,
N\ {\rm divides}\ c\right\}\ ,\\
\Gamma_1(N)&=&\,\left\langle \Gamma(N),\left(\matrix{1&1\cr 0&1}\right)\right\rangle\ .
\eea
For example, $\H/\Gamma_0(2)$ and $\H/\Gamma(2)$ are spheres with 2 and
3 punctures, respectively, while e.g.\  $\H/\Gamma_0(24)$ is a torus with 7
punctures.

The modular functions for the compact surface $\Sigma=\overline{\H/\Gamma}$ are
easy to describe: they are the meromorphic functions $f$ on $\H$, which are
also meromorphic at the cusps $\Q\cup\{\i\infty\}$, and which obey the
symmetry $f({a\tau+b\over c\tau+d})=f(\tau)$ for all $\left(\matrix{a&b\cr
c&d}\right)\in\Gamma$.  The precise definition of `meromorphic at the cusps'
isn't important here. For example, for any group $\Gamma$ obeying
\be
\left(\matrix{1&t\cr 0&1}\right)\in\Gamma\ {\rm iff}\ t\in\Z\ ,\ee
a meromorphic function $f(\tau)$ with symmetry $\Gamma$ will have a Fourier expansion
$\sum_{n=-\infty}^\infty a_nq^n$ for $q=e^{2\pi\i\tau}$ ($q$ is a local
coordinate for $\tau=\i\infty$); then we say $f$ is meromorphic at the cusp
$\i\infty$ iff all but finitely many $a_n$, for $n<0$, are nonzero.

We say a group $\Gamma$ is genus-0 when  $\Sigma=\overline{\H/\Gamma}$ is a
sphere. For these $\Gamma$, there will be a uniformising function
$f_\Gamma(\tau)$ identifying $\Sigma$ with the Riemann sphere
$\C\cup\{\infty\}$. That is, $f_\Gamma$ will be the mother-of-all modular
functions; i.e., it is a modular function for $\Gamma$, and
any other modular function $f(\tau)$ for $\Gamma$ can be
written uniquely as a rational function $poly(f_\Gamma(\tau))/poly(f_\Gamma(\tau))$.
This function $f_\Gamma$ is not quite unique (SL$_2(\R)$ permutes these generating
functions). By contrast, in genus $>0$ two (non-canonical) generating functions will
be needed.

The groups $\Gamma$ we are interested in are genus-0, contain some $\Gamma_0(N)$,
and obey (2.4). We call such $\Gamma$ {\it genus-0 groups of moonshine-type}. Cummins
\cite{Cu2} has classified all of these --- there are precisely 6486 of them.
For these groups (and more generally a group containing a $\Gamma_1(N)$
rather than a $\Gamma_0(N)$) there is a canonical choice of generator $f_\Gamma$:
we can always choose it uniquely so that it has a $q$-expansion of the
form $q^{-1}+\sum_{n=1}^\infty a_nq^n$. This choice of generator is called
the {\it Hauptmodul}, and will be denoted $J_\Gamma(\tau)$. Some examples
are
\bea
J_{\Gamma(1)}(\tau)&=&q^{-1}+196884\,q+21493760\,q^2+864299970\,q^3+\cdots\\
J_{\Gamma_0(2)}(\tau)&=&q^{-1}+276q-2048q^2+11202q^3-49152q^4+184024q^5+\cdots\\
J_{\Gamma_0(13)}(\tau)&=&q^{-1}-q+2\,q^2+q^3+2\,q^4-2\,q^5-2\,q^7-2\,q^8+q^9+\cdots\\
J_{\Gamma_0(25)}(\tau)&=&q^{-1}-q+q^4+q^6-q^{11}-q^{14}+q^{21}+q^{24}-q^{26}+\cdots\ .\eea
Of course $J_{\Gamma(1)}$ is the famous $J$-function.
Exactly 616 of these Hauptmoduls have integer coefficients (171 of which are
the McKay--Thompson series); the remainder have cyclotomic integer coefficients. 

\section{The algebraic meaning of the Hauptmodul property}

The conceptual gap will be bridged only when we can directly relate the
definition of a Hauptmodul (which is inherently topological), with the 
recursions, like (1.2), coming from the twisted denominator identities.

The easiest way to produce functions invariant with respect to some symmetry,
is to average over the group. For example, given any function $f(x)$,
the average $f(x)+f(-x)$ is invariant under $x\leftrightarrow -x$. When
the group is infinite, a little more subtlety is required but the same idea
can work.

For example, take $\Gamma={\rm SL}_2(\Z)$ and let $p$ be any prime. Then
\bea
\Gamma \left(\matrix{p&0\cr 0&1}\right)\Gamma&=&\{A\in M_{2\times 2}(\Z)\,|\,
{\rm det}(A)=p\}\nonumber\\ &=&\Gamma\left(\matrix{p&0\cr 0&1}\right)\cup\bigcup_{k=0}^{p-1}
\Gamma\left(\matrix{1&k\cr 0&p}\right)\ .\eea
This means that, for any modular function $f(\tau)$ of SL$_2(\Z)$,
$f(p\tau)$ will no longer be SL$_2(\Z)$-invariant, but
\be s_f^{(p)}(\tau):=f(p\tau)+\sum_{k=0}^{p-1}f\left({\tau+k\over p}\right)\ee
is. Considering now $f$ to be the Hauptmodul $J$, we thus obtain that $s_J^{(p)}(\tau)=P(J(\tau))/Q(J(\tau))$, for polynomials
$P,Q$. By considering poles and the surjectivity of $J$, we see that $Q$
must be constant, and hence that $s_J^{(p)}$ must be a polynomial in $J$.
The same will hold for any $s_{J^\ell}^{(p)}$.

This implies that there is a monic polynomial $F_p(x,y)$ of degree
$p+1$ in $x,y$, such that
\be 
F_p(J(\tau),J(p\tau))=F_p\left(J(\tau),J\left({\tau+k\over p}\right)\right)=0\ee 
for all $k=0,\ldots,p-1$, or equivalently
\be
F_p(J(\tau),Y)=(J(p\tau)-Y)\prod_{k=0}^{p-1}\left(J\left({\tau+k\over p}\right)-Y\right)\ . \ee
For example,
\bea
F_2(x,y)&=&(x^2-y)(y^2-x)-393768\,(x^2+y^2)-42987520\,xy\nonumber\\
&&-40491318744\,(x+y)
+12098170833256\ .\eea

There is nothing terribly special about $p$ being prime; for a
composite number $m$, the sum in e.g.\ (3.2)
 becomes a sum over $\left(\matrix{m/d&k\cr 0& d}\right)$, for all
divisors $d$ of $m$ and all $0\le k<d$. Write ${\cal A}_m$ for the set
of all these pairs $(d,k)$. Note that its cardinality $\|{\cal A}_m\|$ is 
$\psi(m)=m \prod_{p|m}(1+1/p)$.

\medskip\noindent{\bf Definition 1.} {\it Let $h(\tau)=q^{-1}+\sum_{n=1}^\infty a_nq^n$. We say that $h(\tau)$
 satisfies a modular equation of order $m>1$, if there is a monic polynomial
 $F_m(x,y)\in\C[x,y]$ such that $F_m$ is of degree
 $\psi(m)$ in both $x$ and $y$, and}
\be
F_m(h(\tau),Y)=\prod_{(d,k)\in{\cal A}_m} \left(h\left({m\tau\over d^2}+{k\over
d}\right)-Y\right)\ .\ee

In the following, 
it is unnecessary to assume that the series $h$ converges; it is enough to
require that (3.6) holds formally at the level of $q$-series. An easy consequence
of this definition is that $F_m(x,y)=F_m(y,x)$. 

We've learnt above that
the Hauptmodul $J$ satisfies a modular equation of all orders $m>1$. In fact
similar reasoning verifies, more generally, that:

\medskip\noindent{\bf Proposition 1.} (a) \cite{CuN} {\it If $J_\Gamma(\tau)$ is the Hauptmodul of some genus-0 group
$\Gamma$ of moonshine-type, with rational coefficients, then $J_\Gamma$
satisfies a modular equation for all $m$ coprime to $N$.}

\smallskip\noindent{(b)} {\it Likewise, any
McKay--Thompson series $T_g(\tau)$ satisfies a modular equation for any
$m$ coprime to the order of the element $g\in\M$.}\medskip

Recall there are 616 such $J_\Gamma$, and 171 such $T_g$. The $N$ in part (a)
 is the level of any
congruence group $\Gamma_0(N)$ contained  in $\Gamma$. Part (b) involves
showing that the recursions such as (1.2) imply the modular equation property
$q$-coefficient-wise.

Note that the
Hauptmodul property of $J_\Gamma$ plays a crucial role in the proof that
they satisfy modular equations. Could the converse of the Proposition hold?

Unfortunately, that is too naive. In particular, $h(\tau)=q^{-1}$ also
satisfies a modular equation for any $m>1$: take $F_m(x,y)=(x^m-y)(y^m-x)$.
Using Tchebychev polynomials, it is easy to show  that $h(\tau)=q^{-1}+q$
(which is essentially cosine) likewise satisfies a modular equation for any $m$.

However, Kozlov, in a thesis directed by Meurman, proved the following
remarkable fact:

\medskip\noindent{\bf Theorem 2.} \cite{Ko} {\it If $h(\tau)=q^{-1}+\sum_{n=1}^\infty a_nq^n$ satisfies a modular
equation for all $m>1$, then either $h=J$, $h(\tau)=q^{-1}$, or $h(\tau)=q^{-1}\pm q$.}\medskip

His proof breaks down when we no longer have all those modular equations,
but it gives us confidence to hope that modular equations
could provide an algebraic interpretation of what it means to be a Hauptmodul.
Indeed that is the case!

\medskip\noindent{\bf Theorem 3.} \cite{CG} {\it Suppose a formal series
$h(\tau)=q^{-1}+\sum_{n=1}^\infty a_nq^n$ satisfies a modular
equation for all $m\equiv 1$ (mod $K$). Then $h(\tau)$ is holomorphic throughout
$\H$. Write}
\be \Gamma_h:=\left\{\left(\matrix{a&b\cr
c&d}\right)\in{\rm SL}_2(\R)\,|\,h\left({a\tau+b\over c\tau+d}\right)=h(\tau)\ \forall 
\tau\in\H\right\}\ .\ee 

\noindent{(a)} {\it If $\Gamma_h\ne\left\{\pm\left(\matrix{1&n\cr 0&1}
\right)\,|\,n\in\Z\right\}$, then $h$ is a Hauptmodul for $\Gamma_h$, and $\Gamma_h$
obeys (2.4) and contains $\Gamma_0(N)$ for some $N|K^\infty$}.

\smallskip\noindent{(b)} {\it If $\Gamma_h=\left\{\pm\left(\matrix{1&n\cr 0&1}
\right)\,|\,n\in\Z\right\}$, and the coefficients $a_n$ of $h$ are algebraic integers,
then $h(z)=q^{-1}+\xi q$ where $\xi=0$ or $\xi^{{\rm gcd}(K,24)}=1$.}\medskip

By `$N|K^\infty$' we mean any prime dividing $N$ also divides $K$.
 Of course part (a) implies that $\Gamma_h$ is genus-0 and of moonshine-type.
When $h=T_g$, $K=o(g)$ works (see Prop.1(b)), and all coefficients are
integers, and so Theorem 3 establishes the Hauptmodul property and fills
the conceptual gap. The proof of Theorem 3 is difficult: if $h(\tau_1)=h(\tau_2)$,
then it is fairly easy to prove that locally there is an invertible holomorphic
map $\alpha$ sending an open disc about $\tau_1$ onto one about $\tau_2$; the
hard part of the proof is to show that $\alpha$ extends to a globally
invertible map $\H\rightarrow\H$ (and hence lies in $\Gamma_h$).

The converse of Theorem 3 is also true:

\medskip\noindent{{\bf Proposition 2.}} \cite{CG} {\it If $h(\tau)=q^{-1}+
\sum_{n=1}^\infty a_nq^n$ is a Hauptmodul for a group $\Gamma_h$ of
moonshine-type, and the coefficients $a_n$ all lie in the cyclotomic field
$\Q[\xi_N]$, then there exists a generalised modular equation
for any order $m$ coprime to $N$. Moreover, the field generated over $\Q$
by all coefficients $a_n$ will be a Galois extension of $\Q$, with Galois
group of exponent 2.}\medskip

The exponent 2 condition  means that that field is generated over $\Q$ by
a number of square-roots of rationals. We write $\xi_N$ for the root
of unity $\exp[2\pi\i/N]$. The condition that all $a_n$ lie
in the cyclotomic field should be automatically satisfied. By a `generalised
modular equation' of order $m>1$, we mean that there is a polynomial
 $F_m(x,y)\in\Q[\xi_N][x,y]$ such that $F_m$ is monic of degree
 $\psi(m)$ in both $x$ and $y$, and
\be
F_m((\sigma_m.h)(\tau),Y)=\prod_{(d,k)\in{\cal A}_m} \left(h\left({m\tau\over
d^2}+{k\over d}\right)-Y\right)\ .\ee
We also have the symmetry condition $F_m(x,y)=(\sigma_m.F_m)(y,x)$. Here,
$\sigma_m\in{\rm Gal}(\Q[\xi_N]/\Q)\cong\Z_N^*$ is the Galois
automorphism sending $\xi_N$ to $\xi_N^m$; it acts on
$h$ and $F_m$ coefficient-wise. The beautiful relation of modular functions to 
cyclotomic fields and their Galois groups is classical and is reviewed in
e.g.\ Chapter 6 of  \cite{La}.

Proposition 2 explains why Proposition 1 predicts more modular equations than
Theorem 3 assumes: if all coefficients
$a_n$ in Theorem 3 are rational, then indeed we'd get that $h$ would satisfy an 
ordinary modular equation for all $m$ coprime to $n$.

The lesson of Moonshine is that we probably shouldn't completely ignore
the exceptional functions in Theorem 3(b). It is tempting to call those
25 functions {\it modular fictions} (following John McKay). So a question could be:

\medskip\noindent{{\bf Question 1.}} {\it What is the question in e.g.\
vertex algebras, for which the modular fictions are the answer?}\medskip

Theorem 3 requires many more modular equations than is probably necessary.
 In particular, the computer experiments in \cite{CM} show that if $h$ has
integer coefficients and satisfies modular equations of order 2 and 3, then
$f$ is either a Hauptmodul, or a modular fiction. Cummins  has
made the following conjecture:

\medskip\noindent{{\bf Conjecture 1.}} \cite{Cu1} {\it Let $p_1,p_2$ be
distinct primes, and $a_i\in\C$. Suppose $h(\tau)=q^{-1}+\sum_{n=1}^\infty 
a_nq^n$ satisfies modular equations of order $p_1,p_2$. Then }

\noindent{(a)} {\it If $\Gamma_h\ne\left\{\pm\left(\matrix{1&n\cr 0&1}
\right)\,|\,n\in\Z\right\}$, then $h$ is a Hauptmodul for $\Gamma_h$, and $\Gamma_h$
obeys (2.4) and contains $\Gamma_1(N)$ for some $N$ coprime to $p_1,p_2$}.

\smallskip\noindent{(b)} {\it If $\Gamma_h=\left\{\pm\left(\matrix{1&n\cr 0&1}
\right)\,|\,n\in\Z\right\}$,
then $h(z)=q^{-1}+\xi q$ where $\xi=0$ or $\xi^{{\rm gcd}(p_1-1,p_2-1)}=1$.}
\medskip

We are far from proving this. However, 
if $h$ obeys a modular equation of order $m$ for all $m$ with the property that all
prime divisors $p$ of $m$ obey $p\equiv 1$ (mod $K$) for some fixed $K$, then
$h$ is either a Hauptmodul for $\Gamma_h$ containing some $\Gamma_1(N)$,
or $h$ is `trivial'  (see \cite{Cu1} for details and a proof). The converse
again is known to be true (again provided the $a_i$ are cyclotomic).

It would be interesting to apply similar arguments to fill the related
conceptual gaps of Modular Moonshine \cite{B3} and Baby Moonshine \cite{Ho}.
Modular equations have many uses in number theory, besides these in Moonshine
--- see e.g.\ \cite{Co} for important applications to class field theory.
Modular equations are also closely related to the notion of replicable functions 
(see e.g.\ \cite{Nor}).

\section{The meaning of moonshine}

As mentioned earlier, the greatest open challenge for Monstrous Moonshine is
to find a second independent proof. In this section we briefly sketch some thoughts on
what this proof may involve; see \cite{Ga2} for details.

A powerful guide to Monstrous Moonshine has been rational conformal field theory (RCFT).
Modularity arises in RCFT through the conjunction of two standard pictures:

\begin{itemize}

\item[(1)] {\it canonical quantisation} presents us with a state space $V$, 
carrying a representation of the symmetries of the theory, a Hamiltonian operator
$H$, etc. In RCFT, the quantum amplitudes involve graded traces such as Tr$_Vq^H$, defining the coefficients
of our $q$-expansions.

\item[(2)] The {\it Feynman picture} interprets the amplitudes using path integrals.
In RCFT this permits us to interpret these graded traces as functions (sections) over
moduli spaces, and hence they carry actions by the relevant mapping class groups
such as SL$_2(\Z)$. This gives us modularity.
\end{itemize}

In Monstrous Moonshine, canonical quantisation is successfully abstracted into
the language of vertex operator algebras (VOAs). The present proof of the Conway--Norton
conjectures however ignores the Feynman side, and with it the lesson from RCFT
that modularity is ultimately topological. Perhaps this is where to search for a
second more conceptual proof. After all, the proof of the modularity of VOA
characters \cite{Zh} --- perhaps the deepest result concerning VOAs --- follows exactly
this RCFT intuition. Let's briefly revisit the RCFT treatment of characters.

In an RCFT, with `chiral algebra' (i.e.\ VOA) $V$, the character of a `sector'
(i.e.\ $V$-module)
$M$ is essentially the amplitude associated to a torus with one field (state)
inserted --- we call it a `one-point function' on the torus. Fix a torus $\C/(\Z+\Z \tau)$,
a local parameter $z\in\C$ at the marked point (which we can take to be 0), and the state $v$ which we're inserting at 0 ($v$ can belong to any $V$-module, but for now
we'll take $v\in V$). The local parameter $z$ is needed
for sewing surfaces together at the marked points (a fundamental process in 
RCFT). We get a moduli 
space $\widehat{\cal M}_{1,1}$ of `extended once-marked tori',
i.e.\ tori with a choice of local parameter $z$ at 0. 
In this language the conformal symmetry of the RCFT becomes actions of
the Virasoro algebra; the actions of this infinite-dimensional Lie algebra
on the extended moduli spaces $\widehat{\cal M}_{g,n}$ are responsible for 
much of the special mathematical features of RCFT.

For convenience assume that $Hv=kv$ (this eigenvalue
$k\in\Q$ is called the `conformal weight' of $v$). The character is given by
\be\chi_M(\tau,v,z):={\rm Tr}_M Y(v,e^{2\pi\i z})\,q^{H-c/24}=e^{-2k\pi\i z}
\,{\rm Tr}_Mo(v)\,q^{H-c/24}\ ,\ee
where $c$ is the `central charge' ($c=24$ in Monstrous Moonshine), and $o(v)$
is an endomorphism commuting with $H$ (also called $L_0$). 
What naturally acts on
these $\chi_M$ is the mapping class group $\widehat{\Gamma}_{1,1}$ of
$\widehat{\cal M}_{1,1}$.  

This extended moduli space $\widehat{\cal M}_{1,1}$
is much larger than the usual moduli space  ${\cal M}_{1,1}=\H/{\rm SL}_2(\Z)$ of
a torus with one marked point, and the mapping class group  $\widehat{\Gamma}_{1,1}$
is larger than the familiar mapping class group  ${\Gamma}_{1,1}={\rm SL}_2(\Z)$.
In fact, $\widehat{\Gamma}_{1,1}$ can be naturally identified with the braid
group
\be {\cal B}_3=\langle \sigma_1,\sigma_2\,|\,\sigma_1\sigma_2\sigma_1=
\sigma_2\sigma_1\sigma_2\rangle\ ,\ee
and acts on the characters by
\bea\sigma_1.\chi_M(\tau,v,z)&=&e^{-2\pi\i k/12}\,\chi_M(\tau+1,v,z)\ ,\\
\sigma_2.\chi_M(\tau,v,z)&=&e^{-2\pi\i k/12}\,\chi_M\left({\tau\over 1-\tau},{v\over(1-\tau)^k},
z\right)\ .\eea
Thus in RCFT it is really ${\cal B}_3$ and not SL$_2(\Z)$ which acts on the
characters. This is
usually ignored because we specialise $\chi_M$, and more fundamentally because
typically we consider only insertions $v\in V$, and what results is a true 
action of the modular group SL$_2(\Z)$. But taking $v$ from other
$V$-modules is equally fundamental in the theory, and for those insertions
we only get a projective action of SL$_2(\Z)$ (though again a true action of
${\cal B}_3$). 

This is just a hint of a much more elementary phenomenon. Recall that a
modular form $f$ for $\Gamma:={\rm SL}_2(\Z)$ is a holomorphic function $f:\H\rightarrow
\C$, which is also holomorphic at the cusps, and which obeys
\be f\left(\frac{a\tau+b}{c\tau+d}\right)=\mu\left(\matrix{a&b\cr c&d}\right)\,(c\tau+d)^k
f(\tau)\ \qquad\forall\left(\matrix{a&b\cr c&d}\right)\in\Gamma\ , \ee
for some $k\in\Q$ (called the {\it weight}) and some function $\mu$ (called
the {\it multiplier}) with modulus $|\mu|=1$. For example, the Eisenstein
series
\be E_k(\tau)=\sum_{(m,n)\in\Z^2}\!\!\!{}'\,(m\tau+n)^{-k}\ee
for even $k>2$ is a modular form of weight $k$ with trivial multiplier $\mu$, but the
Dedekind eta
\be \eta(\tau)=q^{1/24}\prod_{n=1}^\infty (1-q^n)\ee
is a modular form of weight $k=1/2$ with a nontrivial multiplier, given by
\be \mu\left(\matrix{a&b\cr c&d}\right)=
\exp\left(\pi\i \,\left({a+d\over 12c}-{1\over 2}- 
\sum_{i=1}^{c-1}{i\over c}\,\left({di\over c}-\left\lfloor {di\over c}\right\rfloor
-{1\over 2}\right)\right)\right)\ \ee 
when $c>0$. 

$\H$ can be regarded as a homogeneous space SL$_2(\R)/{\rm SO}_2(\R)$. Nowadays
we are taught to lift a modular form $f$ from $\H$ to SL$_2(\R)$: 
\be\phi_f\left(\matrix{a&b\cr c&d}\right):=f\left({a\i+b\over c\i+d}\right)\,
(c\i+d)^{-k}\,{\mu}\left(\matrix{a&b\cr c&d}\right)^*\ .\ee
We've sacrificed the implicit SO$_2(\R)$-invariance and explicit 
$\Gamma$-covariance of $f$, for explicit SO$_2(\R)$-covariance
and explicit $\Gamma$-invariance of $\phi_f$. This is significant, because 
compact Lie groups
like the circle SO$_2(\R)$ are much easier to handle than infinite discrete groups 
like SL$_2(\Z)$. The result is a much more conceptual and powerful picture.

Thus a modular form should be regarded as a function on the orbit space
$X:=\Gamma\backslash{\rm SL}_2(\R)$. Remarkably, this 3-space $X$ can be naturally
identified with the complement of the trefoil! We are thus led to ask:

\medskip\noindent{{\bf Question 2.}} {\it Do modular forms for} SL$_2(\Z)$
{\it see the trefoil?}\medskip

An easy calculation shows that the fundamental group $\pi_1(X)$ is in fact
the braid group ${\cal B}_3$! It is a central extension of SL$_2(\Z)$ by
$\Z$. In particular, the quotient of ${\cal B}_3$ by its centre $\langle
(\sigma_1\sigma_2\sigma_1)^2\rangle$ is PSL$_2(\Z)$; the isomorphism
${\cal B}_3/\langle(\sigma_1\sigma_2\sigma_1)^4\rangle\cong {\rm SL}_2(\Z)$
is defined by the (reduced and specialised) Burau representation
\be\sigma_1\mapsto \left(\matrix{1&1\cr 0&1}\right)\ ,\qquad
\sigma_2\mapsto\left(\matrix{1&0\cr -1&1}\right)\ .\ee
Through this map, which is implicit in (4.3) and (4.4), ${\cal B}_3$ acts on modular forms, and the multiplier $\mu$ can be lifted to ${\cal B}_3$. For example,
the multiplier of the Dedekind eta becomes
\be \mu(\beta)=\xi_{24}^{{\rm deg}\,\beta}\qquad \forall\beta\in{\cal B}_3\ ,\ee
where `deg$\,\beta$' denotes the crossing number or degree of a braid.
This is vastly simpler than (4.8)!

In hindsight it isn't so surprising that the multiplier is simpler as a
function of braids than of $2\times 2$ matrices. The multiplier $\mu$ will be
 a true representation of SL$_2(\Z)$ iff the weight $k$ is integral; otherwise
it is only a projective representation. And the standard way to handle
projective representations is to centrally extend. Of course number theorists
know this, but have  preferred using the minimal necessary extension; as
half-integer weights are the most common, they typically only look at
a $\Z_2$-extension of SL$_2(\Z)$ called the metaplectic group Mp$_2(\Z)$.
But unlike ${\cal B}_3$, Mp$_2(\Z)$ isn't much different from the modular
group and the multipliers don't simplify much when lifted to Mp$_2(\Z)$.
At least in the context of modular forms, {\it the braid group
can be regarded as the universal central extension of the modular group,
and the universal symmetry of its modular forms.} 

Topologically, SL$_2(\R)$ is the interior of the solid torus, so its
universal covering group $\widetilde{{\mathrm{SL}}_2(\R)}$ will be the
interior of the solid helix, and a central extension by $\pi_1\cong \Z$ of 
SL$_2(\R)$.  $\widetilde{{\mathrm{SL}}_2(\R)}$ can be realised  \cite{LV}
as the set of all pairs $\left(\left(\matrix{a&b\cr c&d}\right),n\right)$
where $\left(\matrix{a&b\cr c&d}\right)\in {\mathrm{SL}}_2(\R)$ and
$n\equiv 0,1,2,3$ (mod 4) depending on whether $c=0$ and $a>0$, $c<0$,  
$c=0$ and $a<0$, or $c>0$, respectively. The group operation is
$(A,m)(B,n)=(AB,m+n+\tau)$, where $\tau\in\{0,\pm 1\}$ is called the
Maslov index. Just as SL$_2(\Z)$ is the set of all integral points in
SL$_2(\R)$, the braid group ${\cal B}_3$ is the set of all integral points in
$\widetilde{{\mathrm{SL}}_2(\R)}$.

Incidentally, similar
comments apply when SL$_2(\Z)$ is replaced with other discrete groups
--- e.g.\ for $\Gamma(2)$ the relevant central extension is the pure braid
group ${\cal P}_3$. It would be interesting to topologically identify
 the central
extension for all the genus-0 groups $\Gamma_g$ of Monstrous Moonshine.

So far we have only addressed the issue of modularity. A more subtle question
in Moonshine is the relation of the Monster to the genus-0 property. Our
best attempt at answering this is that the Monster is probably the largest
exceptional 6-transposition group \cite{N}. This relates to Norton's generalised
 Moonshine through the notion of quilts (see e.g.\ \cite{Hsu}). The relation of
`6' to genus-0 is that $\H/\Gamma(n)$ is genus-0 iff $n<6$, while $\H/\Gamma(6)$
is `barely' genus 1. The notion of quilts, and indeed the notion of generalised
Moonshine and orbifolds in RCFT, is related to braids through the right action
of ${\cal B}_3$ on any $G\times G$ (for any group $G$) given by
\be (g,h).\sigma_1=(g,gh)\ ,\qquad (g,h).\sigma_2=(gh^{-1},h)\ .\ee

Limited space has forced us to be very sketchy here. For more on all these topics,
see \cite{Ga2} (which you are urged to purchase). We suggest that the braid group ${\cal B}_3$ and related
central extensions may play a
central role in  a new, more conceptual proof of  the Monstrous Moonshine
conjectures.

\bigskip\noindent{\bf Acknowledgements.} This paper was written at the
University of Hamburg, whom I warmly thank for their hospitality. My research
was supported in part by the von Humboldt Foundation and by NSERC.

\end{document}